\documentclass[11pt]{amsart}
\usepackage[T1]{fontenc}

\usepackage[top=2cm,bottom=2cm,left=3cm,right=3cm,marginparwidth=1.75cm]{geometry}

\usepackage[numbers,sort&compress,square]{natbib}
\usepackage[colorlinks=true,citecolor=red,linkcolor=blue,urlcolor=blue]{hyperref}
\usepackage{amsfonts,bm,mathtools,amssymb}
\usepackage{xcolor}

\usepackage{xparse}

\usepackage[shortlabels]{enumitem}

\usepackage{amsthm}
\usepackage[nameinlink,noabbrev]{cleveref}
\usepackage{aliascnt}

\def\eps{\varepsilon}

\newtheorem{theorem}{Theorem}[section]

\newaliascnt{lemma}{theorem}
\newtheorem{lemma}[lemma]{Lemma}
\aliascntresetthe{lemma}
\crefname{lemma}{Lemma}{Lemmas}
\Crefname{lemma}{Lemma}{Lemmas}

\newaliascnt{proposition}{theorem}
\newtheorem{proposition}[proposition]{Proposition}
\aliascntresetthe{proposition}
\crefname{proposition}{Proposition}{Propositions}
\Crefname{proposition}{Proposition}{Propositions}

\newaliascnt{corollary}{theorem}
\newtheorem{corollary}[corollary]{Corollary}
\aliascntresetthe{corollary}
\crefname{corollary}{Corollary}{Corollaries}
\Crefname{corollary}{Corollary}{Corollaries}

\theoremstyle{definition}
\newaliascnt{definition}{theorem}

\aliascntresetthe{definition}
\crefname{definition}{Definition}{Definitions}
\Crefname{definition}{Definition}{Definitions}

\theoremstyle{plain}
\newaliascnt{example}{theorem}
\newtheorem{example}[example]{Example}
\aliascntresetthe{example}
\crefname{example}{Example}{Examples}
\Crefname{example}{Example}{Examples}

\theoremstyle{definition}
\newaliascnt{assumption}{theorem}

\aliascntresetthe{assumption}
\crefname{assumption}{Assumption}{Assumptions}
\Crefname{assumption}{Assumption}{Assumptions}

\theoremstyle{plain}
\newaliascnt{remark}{theorem}
\newtheorem{remark}[remark]{Remark}
\aliascntresetthe{remark}
\crefname{remark}{Remark}{Remarks}
\Crefname{remark}{Remark}{Remarks}

\allowdisplaybreaks[4]

\DeclareMathOperator{\Cov}{Cov}
\DeclareMathOperator{\Var}{Var}

\DeclarePairedDelimiterXPP\pk[1]{\mathbb{P}}\{\}{}{#1}
\DeclarePairedDelimiterXPP\E[1]{\mathbb{E}}\{\}{}{#1}

\newcommand{\R}{\mathbb R}

\newcommand{\M}{\mathcal M}

\providecommand{\supp}{\operatorname{supp}}

\newcommand{\bqn}[1]{\begin{equation}#1\end{equation}}

\def\T{\mathcal T}
\def\PT{p_\T(u)}
\def\IX{\mathfrak m_X}
\def\IZ{\mathfrak m_{Z_*}}

\title[Tail behaviour of Gaussian Infimum]{Tail behaviour of infimum of Gaussian processes}

\author{Krzysztof D\c{e}bicki}
\address{Krzysztof D\c{e}bicki, Mathematical Institute, University of Wroc\l aw, pl. Grunwaldzki 2/4, 50-384 Wroc\l aw, Poland}
\email{Krzysztof.Debicki@math.uni.wroc.pl}

\author{Enkelejd Hashorva}
\address{Enkelejd Hashorva, Department of Actuarial Science, University of Lausanne, UNIL-Dorigny, 1015 Lausanne, Switzerland}
\email{Enkelejd.Hashorva@unil.ch}

\author{Svyatoslav Novikov}
\address{Svyatoslav Novikov, Department of Actuarial Science, University of Lausanne, UNIL-Dorigny, 1015 Lausanne, Switzerland}
\email{Svyatoslav.Novikov@unil.ch}

\date{\today}

\begin{document}

\begin{abstract}
In this contribution we investigate the upper tail of
$\IX=\inf_{t\in\mathcal T}X(t)$,
where $X$ is a centred Gaussian process (GP) with continuous sample paths
on a nonempty compact metric space $\mathcal T$.
As suggested by the findings of \cite{chakrabarty2018asymptotic} the
exponential rate is governed by an optimal-measure problem, while the exact
prefactor depends on the local geometry of the corresponding contact set.
In the full-contact regime we strengthen the known Gumbel
max-domain of attraction property of the distribution of $\IX$ by
establishing the von Mises condition. For fractional Brownian motion on
$[a,b]$, $0<a<b$ with Hurst index $H\in(0,1/2)$ we also determine the order
of the logarithmic correction to the leading Gaussian exponential.
For finite contact sets we distinguish dominated, smooth-critical and
rough-critical regimes and derive exact asymptotics of
$\PT=\pk*{\IX>u}$ as $u\to\infty$ under suitable local assumptions.
We illustrate our findings with examples including fractional Brownian motion
and the fractional Brownian sheet.
\end{abstract}

\maketitle

\section{Introduction}
\label{1}
Let $X(t), t\in\mathcal T$ be a centred GP with
continuous sample paths on a nonempty compact metric space $\mathcal T$ and let
$R(s,t),s,t\in \mathcal T$ be its covariance function.  Motivated by \cite{chakrabarty2018asymptotic}, we study the extreme value properties of the infimum $\IX =\inf_{t\in\T}X(t)$.
For every nonempty compact $A\subseteq\mathcal T$ write
$p_A(u)=\pk*{\inf_{t\in A}X(t)>u}$.
In order to exclude degenerate cases, we assume throughout that there exists
$u_0>0$ such that for all $u>u_0$
\begin{equation} \label{2}
        \PT  >0.
\end{equation}

We first consider the finite case $\T=\{t_1,\ldots,t_n\}, n\ge 1$, with
 $\Sigma=(R(t_i,t_j))_{1\le i,j\le n}$  nonsingular.  It is well-known that the asymptotics of $p_\T(u)$ as $u\to \infty$  are governed by the quadratic programming problem
\[
        \min_{\bm p\in\Delta_n}\bm p^\top\Sigma\bm p,
        \qquad
        \Delta_n:=\left\{\bm p\in[0,\infty)^n:
        \sum_{i=1}^n p_i=1\right\}
\]
see e.g., \cite{HashorvaHuesler2003}.
Let $\bm p_*$ be the minimiser and put
\[
        \sigma_*^2:=(\bm p_*)^\top\Sigma\bm p_*>0,
        \qquad
        m(t_i):=\frac{(\Sigma\bm p_*)_i}{\sigma_*^2},
        \quad i=1,\ldots,n.
\]
The optimality conditions imply that
\[
        m(t_i)\ge1,\qquad i=1,\ldots,n
\]
with equality whenever the $i$th component of $\bm p_*$ satisfies $p_{*,i}>0$.  Thus if
\[
        S:=\{t_i\in \T:p_{*,i}>0\}
\]
is the active set, $E=\{t_i\in\T:m(t_i)=1\}$ is the contact set, and
$P=E\setminus S$ is the passive set with $\ell$ elements, then
\bqn{\label{3}
        \PT\sim p_E(u)
        \sim C\pk*{\min_{t\in S}X(t)>u},
        \qquad u\to\infty,
}
where $C=1$ if $\ell=0$, while for $\ell>0$, $C$ is the
orthant probability of the centred Gaussian regression residual at the passive
points $t_i\in P$; see \Cref{54} and
\cite[Corollary~4.1]{HashorvaHuesler2003} and
\cite[Example~6.1 and Theorem~5.1]{Hashorva2007}. Moreover,
\cite[Theorem~2.1]{hashorvaNovikov2026highminima} implies that the distribution function $F$ of
$\IX$ belongs to the Gumbel max-domain of attraction (GMDA) with scaling
function $a$ satisfying $a(u)\sim \sigma_*^2/u$ as $u\to\infty$, i.e.,
\begin{equation}
\label{4} \lim_{u\to \infty}
\frac{1-F(u+x a(u))}{1-F(u)} = \exp\left(-x\right), \quad \forall x\in \R.
\end{equation}

For $\T$ a compact interval of $\R$, the GMDA of $F$ was established under
smoothness assumptions in \cite{chakrabarty2018asymptotic}.  For continuous GPs
on general compact metric spaces, the GMDA property under the standing
non-degeneracy assumption \eqref{2} follows from
\cite[Theorem~2.1]{hashorvaNovikov2026highminima}.

When $F$ has an eventually positive differentiable density $F'$,
define its reciprocal hazard by $a(u):=(1-F(u))/F'(u)$ for all sufficiently
large $u$. The von Mises condition
\begin{equation}
\label{5} \lim_{u\to \infty}  a'(u)=0
\end{equation}
does not follow from GMDA and absolute
continuity alone.

For a general compact parameter space $\T$, our first question is:\\

 \textbf{Q1.} \textit{Under what additional conditions does the reciprocal
hazard of $F$ satisfy the von Mises condition \eqref{5}?}\\ 

Next, let $\M_1(\mathcal T)$ denote the
set of Borel probability measures on $\mathcal T$.  For every
$\mu\in\M_1(\mathcal T)$, averaging gives the upper bound
\bqn{\label{6}
        \PT \le \pk*{X_\mu>u},
        \qquad
        X_\mu:=\int_{\mathcal T} X(t) \mu(dt).
}
By weak compactness of $\M_1(\mathcal T)$ and continuity of the covariance function
of $X$, the energy functional $\mu\mapsto\Var(X_\mu)$ is
continuous under weak convergence of probability measures
and the minimum is attained.  Any minimiser $\mu_*\in\M_1(\mathcal T)$ is called
an \emph{optimal measure} and we set
\[
        \sigma_*^2
        :=\min_{\mu\in\M_1(\mathcal T)}\Var(X_\mu).
\]

The non-degeneracy assumption \eqref{2} implies
\begin{equation}
\label{7}
        \sigma_*^2
        =\Var(X_{\mu_*})
        >0.
\end{equation}
Indeed, if $\sigma_*^2=0$, then $X_{\mu_*}=0$ almost surely, which is
incompatible with $ \PT >0$ for $u>u_0$.

Next, let $\mu_*$ be an optimal measure.  In view of the optimality criterion in
\Cref{57} below, its normalised potential
\[
        m(t)
        :=
        \frac{\Cov(X(t),X_{\mu_*})}{\Var(X_{\mu_*})},
        \qquad t\in\mathcal T,
\]
satisfies
\[
        m(t)\ge1,\qquad t\in\mathcal T, \qquad\text{and}\qquad
        m(t)=1
        \qquad\text{on }\supp(\mu_*)=:S
\]
and hence the support $S$, which may depend on the choice of optimal measure,
is contained in the contact set, i.e.,
\[
     S \subset    E:=\{t\in\mathcal T:m(t)=1\}.
\]
The normalised potential $m$ and hence $E$ do not depend on the choice of
optimal measure.  Indeed, if $\mu$ and $\nu$ are optimal, then
\[
 \Var\left(\frac{X_\mu+X_\nu}{2}\right)
 =\sigma_*^2-\frac14\Var(X_\mu-X_\nu).
\]
Since $(\mu+\nu)/2$ is a probability measure, optimality forces
$X_\mu=X_\nu$ almost surely, so their potentials coincide.
Both $E$ and $S$ are compact. Gaussian regression gives
\bqn{\label{8}
        X(t)=m(t)X_{\mu_*}+Z_*(t),
        \qquad t\in\mathcal T,
}
where the residual process $Z_*$ is a centred GP independent of
$X_{\mu_*}\sim N(0,\sigma_*^2)$,
see e.g., \cite{chakrabarty2018asymptotic,hashorvaNovikov2026highminima}.

Thus $S$, $E$, and $\mu_*$ are the continuous counterparts of the active
set, contact set, and optimal weights in the finite-dimensional problem introduced above.

The logarithmic asymptotic of $\PT$ follows directly from \eqref{8}.
Indeed, put $Y:=X_{\mu_*}$. Independence and $m\ge1$ then give,
for every fixed $c>0$ and $u>0$,
\[
 \pk*{Y>u+c}
 \pk*{\|Z_*\|_\infty<c}
 \le p_\T(u)\le p_S(u)\le\pk*{Y>u}.
\]
The support of a centred Gaussian measure on the separable Banach space
$C(\T)$ is a closed linear subspace; see
\cite[Theorem~3.6.1]{Bogachev1998}. Hence $\pk*{\|Z_*\|_\infty<c}>0$.
Thus $\sigma_*^2>0$ also implies \eqref{2} and Gaussian tail asymptotics yield
\bqn{\label{9}
        \lim_{u\to\infty}\frac1{u^2}\log \PT
        =\lim_{u\to\infty}\frac1{u^2}\log p_{S}(u)
        =\lim_{u\to\infty}\frac1{u^2}\log \pk*{X_{\mu_*}>u}
        =-\frac1{2\sigma_*^2}<0.
}

In this Gaussian setting, log-concavity of $p_\T$ together with
\eqref{9} yields the GMDA property; see
\cite[Section~3.2]{hashorvaNovikov2026highminima}.
When $E=\T$ hereafter called the \emph{full-contact regime},
\eqref{8} reduces to
\bqn{\label{10}
        X(t)=X_{\mu_*}+Z_*(t),\qquad t\in\mathcal T.
}
This class of GPs with full contact is broad.  Indeed if $Z_*$ is a centred
GP with continuous sample paths and if
$\mu_*\in\M_1(\T)$ satisfies \eqref{11} almost surely,
then for an independent Gaussian random variable (rv)
$Y\sim N(0,\sigma_*^2)$
 the process
\[
        X(t):=Y+Z_*(t),\qquad t\in\T
\]
has $\mu_*$ as an optimal measure and has full contact set $E=\T$.
In the full-contact regime, integrating \eqref{10} against
$\mu_*$ yields almost surely
\begin{equation}
\label{11}
        \int_\T Z_*(t) \mu_*(dt)=0
\end{equation}
and hence
\[
       \IZ := \inf_{t\in\T}Z_*(t)\le0.
\]
Moreover, we have
\[
        \IX=X_{\mu_*}+\IZ
\]
and consequently as shown in \cite{wu2019high}, for every $u\in\mathbb R$
\bqn{\label{12}
  \PT
  =
  \frac{\exp\left(-u^2/(2\sigma_*^2)\right)}
       {\sqrt{2\pi}\sigma_*}
  \int_0^\infty
  \exp\left(
    -\frac{u\varepsilon}{\sigma_*^2}
    -\frac{\varepsilon^2}{2\sigma_*^2}
  \right)
  \pk*{\IZ> -\varepsilon} d\varepsilon.
}
An immediate Abelian consequence of the above identity is
\[
\lim_{u\to\infty}
\frac{\PT}{\pk*{X_{\mu_*}>u}}
=
\pk*{\IZ = 0}.
\]
If $\pk*{\IZ = 0}$ is positive, then both  the exact asymptotics of $\PT$ and the von Mises condition \eqref{5} hold.

We shall show that \eqref{5} remains valid when
$\pk*{\IZ=0}=0$. To obtain exact tail asymptotics in this case,
\eqref{12} reduces the problem to the behaviour of
$r(\varepsilon):=\pk*{\IZ>-\varepsilon}$ as $\varepsilon\downarrow0$.
Indeed, if, for some ultimately positive function $\mathcal L$ slowly varying at infinity,
\begin{equation}
\label{13}
        r(\eps)\sim K\eps^\alpha \mathcal L(1/\eps),
        \qquad K>0,\quad \alpha\ge 0
\end{equation}
holds as $\eps\downarrow0$, then by Karamata's Abelian theorem for Laplace transforms and \eqref{12} we obtain
\begin{equation}
\label{14}
        \PT
        \sim
        \frac{K\Gamma(\alpha+1)\sigma_*^{2\alpha+1}}
             {\sqrt{2\pi}}
        u^{-(\alpha+1)}\mathcal L(u)
        \exp\left(-\frac{u^2}{2\sigma_*^2}\right),
        \qquad u\to\infty.
\end{equation}
Explicit asymptotics for $r(\varepsilon)$ as
$\varepsilon\downarrow0$ are available in some particular cases; see
\Cref{22} below.

Exact asymptotics of $\PT$ are available in particular in the smooth
settings of \cite{chakrabarty2018asymptotic} and in special Brownian cases; see
\cite{dkebicki2016parisian,debicki2026distributional}.  This leads to our second
question:\\

\textbf{Q2.} \textit{Under what conditions can we derive exact tail
asymptotics of $\PT$ as $u\to\infty$ in finite-contact regimes?}\\

To address Q2, we consider the case where the contact set $E$ is finite.
Every optimal measure $\mu_*$ then has finite support and can be written as
\bqn{\label{15}
        \mu_*=\sum_{i=1}^k p_{*,i}\delta_{t_i},
        \qquad p_{*,i}>0,\qquad i=1,\ldots,k, \qquad \sum_{i=1}^k p_{*,i}=1,
}
where   $\delta_{t_i}$ denotes the Dirac probability measure at $t_i$. The support of $\mu_*$ is $S=\{t_1,\ldots,t_k\}$,  the weight vector is
$\bm p_*=(p_{*,1},\ldots,p_{*,k})^\top$ and we set
$\Sigma_*=(R(t_i,t_j))_{1\le i,j\le k}, \bm X_S=(X(t_1),\ldots,X(t_k))^\top$.  The optimality criterion in
\Cref{57} implies
\[
        \Sigma_*\bm p_*=\sigma_*^2\bm1,
        \qquad
        \bm1=(1,\ldots,1)^\top\in\R^k.
\]
If $\Sigma_*$ is nonsingular, Gaussian regression therefore implies
\bqn{\label{16}
     m(t)=\bm a(t)^\top\bm1\ge1, \qquad X(t)=\bm a(t)^\top\bm X_S+Z(t),
        \qquad
        \bm a(t):=\Sigma_*^{-1}
        (R(t_1,t),\ldots,R(t_k,t))^\top,
        \quad t\in\T,
}
where $Z$ is a centred GP independent of $\bm X_S$.  The support points of the optimal measure are special since
\[
m(t_i)= 1, \qquad Z(t_i)=0, \qquad 1 \le i\le k.
\]
The support points in $S$ are called \textit{active} contacts, while the points
in $P=E\setminus S$ are called \textit{passive} contacts.

In general the two residual processes $Z$ and $Z_*$ differ since
\[
        Z_*(t)-
        Z(t)
        =
        \bm a(t)^\top\bigl(\bm X_S-\bm1 X_{\mu_*}\bigr),
        \qquad t\in\T,
\]
however if $S$ is a singleton, then
\[
        Z(t)=Z_*(t),\qquad t\in\T.
\]

The reduction in \eqref{3} to a positive constant times $p_S(u)$ need not
hold for an infinite parameter set, even when $E$ is finite.  To see this, let $\xi,\eta$ be
independent standard normal rvs and for $\alpha,\lambda>0$ set
\[
        X(t):=\xi\bigl(1+\lambda |t|^\alpha\bigr)+\eta t,
        \qquad t\in \T=[-1,1].
\]
For this process, we have
\[
        \mu_*=\delta_0,
        \qquad
        S=E=\{0\},
        \qquad
        \sigma_*^2=1,
        \qquad
        m(t)-1=\lambda |t|^\alpha,
        \qquad
        Z(t)=\eta t,
        \qquad t\in[-1,1]
\]
and the residual GP $Z$ has mean-square H\"older exponent $\beta=1$.
Direct calculation yields
\[
\lim_{u\to\infty}\frac{\PT}{\pk*{\xi>u}}
=
\begin{cases}
1, & 0<\alpha<2,\\[1mm]
\left(1+\dfrac{1}{2\lambda}\right)^{-1/2}, & \alpha=2,\\[2mm]
0, & \alpha>2.
\end{cases}
\]
The three regimes in the above limit reflect the competition between the rate at which
$m(t)-1$ vanishes near $0$, controlled by $\alpha$, and the residual
smoothness exponent $\beta=1$.

The cases $\alpha<2$ and $\alpha=2$ are covered by
\Cref{27} and the interior version in
\Cref{41}, respectively.  The regime in which
the ratio vanishes, illustrated by $\alpha>2$, is outside the scope of
our general exact-asymptotic results.\\

Critical contacts balance potential separation of order $|h|^{2\beta}$
with residual fluctuations on the scale $|h|^\beta$.  The
\emph{smooth-critical} case has quadratic separation and a linear residual
tangent, as when $\alpha=2$ above.  Critical contacts with
$\beta\in(0,1)$ are called \textit{rough-critical}.  The Brownian endpoint model in
\Cref{42} has $\beta=1/2$ and yields an explicit
nontrivial correction factor.

Passive orthant factors and smooth quadratic corrections already appear under
stronger smoothness assumptions in
\cite[Theorem~4.1 and equation~(4.3)]{chakrabarty2018asymptotic}.
Our finite-contact results use local assumptions and also treat rough tangents.

A further motivation comes from \cite{debicki2026expected}, where
the expected-infimum and persistence asymptotics of Brown--Resnick
stationary log-normal processes exhibit endpoint domination for fractional
Brownian motion with $H>1/2$, in contrast to the Brownian boundary case.
Those results concern drifted processes on growing intervals.  Here we
study high minima on fixed compact sets and explain the corresponding
distinction through the geometry of the contact set.\\ 

Our main contributions are as follows:
\begin{enumerate}
  \item We establish the von Mises condition in the full-contact regime
  and determine the order of the logarithmic correction for fractional
  Brownian motion with $H<1/2$.
  \item We obtain an exact finite-contact reduction for dominated support
  contacts, including the Gaussian orthant factor generated by finitely many
  passive contacts. The result also applies to Gaussian random fields,
  including the fractional Brownian sheet.
  \item For a single endpoint contact, we identify the exact corrections in
  rough-critical and smooth-critical regimes through a tangent process or a
  residual derivative.
\end{enumerate}

The validity of the von Mises condition under the general finite-contact
assumptions is not addressed here.\\ 

The paper is organised as follows followed by a section dealing with the
finite-contact case and presenting several examples.   The proofs are relegated  in
\Cref{44}, while the Appendix reviews  the auxiliary optimal measure
results used in the main text.

\section{Full-Contact Regime}
\label{17}

We retain the notation of the Introduction and assume throughout this section
that $\sigma_*^2>0$ and $E=\T$. The Laplace representation
\eqref{12} connects the tail asymptotics of $\PT$ with small deviations
of the residual process $Z_*$. We first establish the von Mises condition,
then turn to logarithmic refinements and an exact Brownian example.\\ 

Throughout the paper $B_H$, $H\in(0,1)$ denotes a standard fractional
Brownian motion (fBM).

The following result answers Q1 in the full-contact regime.

\begin{theorem}
\label{18}
If $E=\T$ and $\sigma_*^2>0$, then the infimum
$\IX=\inf_{t\in\T}X(t)$
has a strictly positive $C^\infty$-density $f_\T$ on $\mathbb R$ and satisfies
\begin{equation}
\label{19}
               h_\T(u):=\frac{f_\T(u)}{p_\T(u)}
        =
        \frac{u}{\sigma_*^2}+o(1),
        \qquad
        h_\T'(u)
        =
        \frac1{\sigma_*^2}+o(1), \quad u\to \infty
\end{equation}
and
\begin{equation}
\label{20}
           a_\T(u):=\frac{p_\T(u)}{f_\T(u)}
        =
        \frac{\sigma_*^2}{u}+o(u^{-2}),
        \qquad
        a_\T'(u)
        \sim
        -\frac{\sigma_*^2}{u^2}, \qquad u\to\infty.
\end{equation}
\end{theorem}

The logarithmic asymptotic in \eqref{9} was refined for Gaussian
Markov processes in \cite[Theorem~3.1]{wu2019high}.  For the rough fBM considered
below, \cite[Theorem~3.2]{wu2019high} already yields the lower correction bound
once full support of the optimal measure is known.  The matching upper bound
follows from the optimal-measure density bound and an $L^1$
small-deviation estimate.\\ 

For the next proposition, let
\[
\T=[a,b], \qquad 0< a< b< \infty.
\]  If $X=B_H$ with $H\in(0,1/2)$,
\cite[Example~2.7 and Lemma~A.2]{hashorvaNovikov2026highminima} yields
$S=E=[a,b]$ and $\mu_*(dt)\ge c_\mu dt$ on $(a,b)$ for some $c_\mu>0$.
Our next result identifies the order of the correction to \eqref{9}.

\begin{proposition}
\label{21} If
$X=B_H, H\in(0,1/2)$,
then there exist constants $D_\pm=D_\pm(a,b,H)$, with
$0<D_-\le D_+<\infty$ such that for all sufficiently large $u$
\[
        -D_+u^{1/(1+H)}
        \le
        \log\pk*{\inf_{t\in[a,b]}X(t)>u}
        +\frac{u^2}{2\sigma_*^2}
        \le
        -D_-u^{1/(1+H)}.
\]
\end{proposition}
 
The following example treats two intervals and makes the dependence of
the prefactor on their lengths and separation explicit.

\begin{example}
\label{22}
Let $B:=B_{1/2}$ be a standard Brownian motion, fix
$0<a_1<b_1<a_2<b_2<\infty$ and set
\[
        \mathcal K:=[a_1,b_1]\cup[a_2,b_2],
        \qquad
        p_{\mathcal K}(u):=
        \pk*{\inf_{t\in\mathcal K}B(t)>u}.
\]
The normalised potential of $\delta_{a_1}$ is $m\equiv1$ on
$\mathcal K$, so $E=\mathcal K$.
The measure $\delta_{a_1}$ is nevertheless the unique optimal measure,
because for every $\mu\in\mathcal M_1(\mathcal K)$
\[
        \Var\left(\int_{\mathcal K}B(t) \mu(dt)\right)
        =a_1+\int_{a_1}^{b_2}\mu([x,b_2])^2 dx
        \ge a_1,
\]
with equality only for $\mu=\delta_{a_1}$ and hence
\[
        \sigma_*^2=a_1,
        \qquad
        Z_*(t)=B(t)-B(a_1).
\]

Set
\[
        A:=b_1-a_1,
        \qquad G:=a_2-b_1,
        \qquad C:=b_2-a_2,
\]
and let
\[
        W(t):=B(a_1+t)-B(a_1),\quad
        \varphi_t(x):=
        \frac{1}{\sqrt{2\pi t}}\exp\left(-\frac{x^2}{2t}\right),
        \qquad
        \Phi(x):=\int_{-\infty}^x\varphi_1(y) dy.
\]
By the reflection principle   and the Markov property  of the Brownian motion $W$
\[
\begin{aligned}
        r_{\mathcal K}(\eps)
        &:={}
        \pk*{
        \inf_{t\in[0,A]\cup[A+G,A+G+C]}W(t)>-\eps
        } \\
        &={}
        \int_0^\infty
        (\varphi_A(y-\eps)-\varphi_A(y+\eps))
        \int_0^\infty
        \varphi_G(z-y)
        \left\{2\Phi\left(\frac{z}{\sqrt C}\right)-1\right\}
         dz dy.
\end{aligned}
\]
Consequently, by dominated convergence
$
        r_{\mathcal K}(\eps)
        \sim \kappa_{A,G,C}\eps,
         \eps\downarrow0,$
where
\begin{align}
                \kappa_{A,G,C}
        &=
        \int_0^\infty\frac{2y}{A}\varphi_A(y)
        \int_0^\infty
        \varphi_G(z-y)
        \left\{2\Phi\left(\frac{z}{\sqrt C}\right)-1\right\}
         dz dy\notag \\
        &= \label{23}
 \frac{2\varphi_A(0)}{\pi}\arctan\sqrt{\frac GC}
  +\varphi_T(0)\left(1+\frac2\pi\arcsin\rho\right)
        \in(0,\infty),
        \end{align}
where $T:=A+G+C$ and $\rho:=\sqrt{AC/((A+G)(G+C))}$.
 
\end{example}

\section{Finite Contact Sets}
\label{24}

We now consider $\T=K$, where $K$ is a nonempty compact subset of
$\mathbb R^d$, $d\ge1$ and assume that the contact set is finite with
\[
        S=\{t_1,\ldots,t_k\},
        \qquad
        E=\{t\in K:m(t)=1\}=S\cup P,
        \qquad
        P=\{s_1,\ldots,s_\ell\}.
\]
We shall suppose that $\mu_*$ is the unique optimal probability measure.
Write $\mu_*$ as in \eqref{15}.
The finite-support optimality relations yield
\[
        \Sigma_*\bm p_*=\sigma_*^2\bm1
\]
and $\Sigma_*$ is nonsingular.  Indeed, if $\Sigma_*\bm v=0$, then
\[
        \bm1^\top\bm v
        =\sigma_*^{-2}\bm p_*^\top\Sigma_*\bm v=0
\]
and $\bm p_*+\varepsilon\bm v$ is a distinct optimal probability vector for
all sufficiently small nonzero $\varepsilon$, unless $\bm v=0$.
Consequently, the regression representation \eqref{16} applies. At the support points
\[
        \bm a(t_i)=\bm e_i,
        \qquad
        Z(t_i)=0,
        \qquad i=1,\ldots,k,
\]
where $\bm e_i$ is the $i$-th standard basis vector.  Setting
\[
        \bm S_u:=u(\bm X_S-u\bm1)
\]
yields for every $u>0$
\bqn{\label{25}
        X(t)-u
        =
        u(m(t)-1)
        +Z(t)
        +\frac1u\bm a(t)^\top\bm S_u, \qquad t\in \T.
}
Thus points outside $E$ have a positive deterministic separation
$m(t)-1$. At a passive contact $s\in P$ by contrast $m(s)=1$ but
$Z(s)$ need not vanish, whereas $Z(t_i)=0$ at every support contact.

\subsection{Dominated contacts}
A support point $t_i$ is called \emph{dominated} if there exist finite
positive constants $\beta_i,C_i,\lambda_i$ and $0<q_i<2\beta_i$ such that
for all sufficiently small admissible $h\in\mathbb R^d$
\begin{equation}
\label{26}
  m(t_i+h)-1\ge \lambda_i\|h\|^{q_i},
\end{equation}
and for all sufficiently small admissible $h,h'\in\mathbb R^d$
\[
        \Var(Z(t_i+h)-Z(t_i+h'))
        \le C_i\|h-h'\|^{2\beta_i}.
\]
Here admissible means that every displayed point belongs to $K$.

\begin{theorem}
\label{27}
Under the standing assumptions of this section if every
$t_i\in S$ is dominated, then as $u\to\infty$
\[
\begin{aligned}
        \pk*{\inf_{t\in K}X(t)>u}
        &\sim C \pk*{X(t_i)>u,\ i=1,\ldots,k} \sim
        \frac{C \prod_{i=1}^k\theta_i^{-1}}
             {(2\pi)^{k/2}\sqrt{\det\Sigma_*}}
        u^{-k}\exp\left(-\frac{u^2}{2\sigma_*^2}\right),
\end{aligned}
\]
where $\theta_i=p_{*,i}/\sigma_*^2$, with the convention $C=1$ when
$P=\varnothing$ and otherwise
\[
        C:=\pk*{Z(s)>0,\ \forall s\in P}>0.
\]
\end{theorem}

\begin{remark}
\label{28}
The regression decomposition and independence yield for $s,t\in\T$
\[
        \Var\bigl(Z(t)-Z(s)\bigr)
        =
        \Var\bigl(X(t)-X(s)\bigr)
        -
        \Var\left(
        \bigl(\bm a(t)-\bm a(s)\bigr)^\top\bm X_S
        \right)
        \le
        \Var\bigl(X(t)-X(s)\bigr).
\]
Consequently, in the definition of a dominated contact it is enough to verify
the corresponding local mean-square H\"older bound for $X$.
\end{remark}

\begin{remark}
\label{29}
Let $\T=[a,b]$, where $0<a<b$, and suppose that $R(a,a)>0$.  If
\[
        R(t,a)>R(a,a),
        \qquad t\in(a,b],
\]
then $\delta_a$ is the unique optimal measure by
\Cref{58} and
\[
        S=E=\{a\},
        \qquad
        m(t)=\frac{R(t,a)}{R(a,a)}.
\]

If $X$ has stationary increments, $X(0)=0$ almost surely, and
$V(t):=\Var(X(t))$ the preceding condition is equivalent to
\[
        V(a+h)>V(a)+V(h),
        \qquad 0<h\le b-a.
\]
If the preceding strict inequality holds for every $0<h\le b-a$, and for all sufficiently small $h>0$
\[
        V(a+h)-V(a)-V(h)\ge c_0h^q,
        \qquad
        V(h)\le Ch^{2\beta},
        \qquad q<2\beta,
\]
then $a$ is dominated and
\[
        \pk*{\inf_{t\in[a,b]}X(t)>u}
        \sim
        \pk*{X(a)>u},
        \qquad u\to\infty.
\]
This includes fractional Brownian motion with $H>1/2$, for which
$q=1$ and $\beta=H$.
\end{remark}

We can now compare \Cref{22} with the case $H\in(1/2,1)$.
Consider $B_H$ on the same set
$\mathcal K=[a_1,b_1]\cup[a_2,b_2]$.
The strict covariance inequality in \Cref{29} holds for every $t>a_1$,
so the unique optimal measure on $\mathcal K$ is $\delta_{a_1}$, with
$S=E=\{a_1\}$. Applying \Cref{29} on $[a_1,b_2]$ and bounding
the infimum over $\mathcal K$ between that over $[a_1,b_2]$ and
$B_H(a_1)$ yields
\[
 \pk*{\inf_{t\in\mathcal K}B_H(t)>u}
 \sim\pk*{B_H(a_1)>u}
 \sim\frac{a_1^H}{\sqrt{2\pi} u}
       \exp\left(-\frac{u^2}{2a_1^{2H}}\right),
 \qquad u\to\infty.
\]
Thus the other endpoints and the gap do not enter the leading term,
whereas in the Brownian case of \Cref{22} they affect the prefactor
through $\kappa_{A,G,C}$.

For fBM on a fixed interval $[a,b]$, $0<a<b<\infty$,
\Cref{22,21,29} distinguish three behaviours. For $H<1/2$, the logarithm of the
probability differs from $-u^2/(2\sigma_*^2)$ by a negative term of order
$u^{1/(1+H)}$.  For $H=1/2$ and $H>1/2$, the exact prefactors multiplying
the leading Gaussian exponential are proportional to $u^{-2}$ and
$u^{-1}$, respectively.

We next apply \Cref{27} to several GPs.

\begin{example} 
\label{30}
Fix $0<a<b<\infty$ and let $\T=[a,b]$.  The criterion in
\Cref{29} applies to
\[
        I(t):=\int_0^tY(s) ds,
\]
where $Y$ is a centred stationary Gaussian process with continuous
nonnegative covariance $\rho$ satisfying $\rho(0)>0$, with the integral understood
in mean square.  The process $I$ has stationary increments, and
$\Cov(I(a),I(a+h)-I(a))>0$ for every $h>0$.  Moreover, we have 
\[
 \lim_{h\downarrow0}\frac{\Cov(I(a),I(a+h)-I(a))}{h}
 =\int_0^a \rho(v) dv>0,
 \qquad
 \Var(I(t)-I(s))\le \rho(0)|t-s|^2.
\]
Thus $q=\beta=1$ and \Cref{27} yields
\[
 \pk*{\inf_{t\in[a,b]}I(t)>u}
 \sim\pk*{I(a)>u},\qquad u\to\infty.
\]
\end{example}

\begin{example} 
\label{31}
Fix $0<a<b<\infty$ and let $\T=[a,b]$.  The covariance criterion in
\Cref{29} also applies to
the Riemann--Liouville process
\[
        X_H^{\rm RL}(t)
        :=\int_0^t(t-s)^{H-\frac12} dB(s),
        \qquad H\in(1/2,1),
\]
where $B$ is a standard Brownian motion.  The covariance of
$X_H^{\rm RL}$ is strictly increasing in $t\ge a$ with the first argument
fixed at $a$.  Differentiating the covariance integral and applying
It\^o's isometry yield
\[
 m(a+h)-1\sim\frac{H}{a}h,\quad h\downarrow0,
 \qquad
 \Var(X_H^{\rm RL}(t)-X_H^{\rm RL}(s))
 \le C_H|t-s|^{2H}
\]
for some $C_H<\infty$.  For the increment bound, rescaling the integrals
in It\^o's isometry uses
\[
 \int_0^\infty
 \bigl((v+1)^{H-1/2}-v^{H-1/2}\bigr)^2 dv<\infty,
\]
hence $q=1$ and $\beta=H$ and thus  \Cref{27} yields
\[
 \pk*{\inf_{t\in[a,b]}X_H^{\rm RL}(t)>u}
 \sim\pk*{X_H^{\rm RL}(a)>u},\qquad u\to\infty.
\]
\end{example}

\begin{example} 
\label{32}
Let
\[
        I_H(t):=\int_0^t B_H(r) dr,
        \qquad H\in(0,1),
\]
where $B_H$ is a standard fractional Brownian motion.  Fix
$0<a<b<\infty$ and let $\T=[a,b]$.  Write
\[
        R_H(r,s):=\Cov(B_H(r),B_H(s)).
\]
For $t>a$ we have
\[
\Cov(I_H(a),I_H(t))-\Var(I_H(a))
=
\int_0^a\int_a^t R_H(r,s) ds dr>0
\]
implying that   $\delta_a$ is the unique optimal measure and
\[
        S=E=\{a\},
        \qquad
        P=\varnothing.
\]
Moreover we have
\[
        \sigma_*^2=\frac{a^{2H+2}}{2H+2},
        \qquad
        m(a+h)-1\sim\frac{H+1}{a}h,
        \qquad h\downarrow0,
\]
while for $s,t\in[a,b]$
\[
        \Var(I_H(t)-I_H(s))\le C|t-s|^2.
\]
Hence $q=\beta=1$ implying that  $a$ is dominated and
\[
\begin{aligned}
\pk*{\inf_{t\in[a,b]}I_H(t)>u}
&\sim \pk*{I_H(a)>u} \sim
\frac{a^{H+1}}{2\sqrt{\pi(H+1)} u}
\exp\left(-\frac{H+1}{a^{2H+2}}u^2\right),
\qquad u\to\infty.
\end{aligned}
\]

\end{example}

\begin{example} 
\label{33}
Let $X(t),t\ge 0$ be the centred Gaussian process with covariance
\[
        R(s,t)=\frac{2st}{s+t},
        \qquad s,t>0;
\]
see \cite{LiS01}.  On $\T=[a,b]$, where $0<a<b<\infty$
\[
        m(t)=\frac{R(t,a)}{R(a,a)}
        =\frac{2t}{a+t}>1,
        \qquad t>a.
\]
Thus $\delta_a$ is the unique optimal measure and
\[
        \mu_*=\delta_a,
        \qquad
        S=E=\{a\},
        \qquad
        \sigma_*^2=a.
\]
For the regression residual $Z(t)=X(t)-m(t)X(a)$
\[
        m(a+h)-1=\frac{h}{2a+h}\ge\frac{h}{2b},
        \qquad
        \Var(Z(t)-Z(s))
        \le\frac{(t-s)^2}{t+s}
        \le\frac{(t-s)^2}{2a},
\]
hence $q=\beta=1$ and
\[
        \pk*{\inf_{t\in[a,b]}X(t)>u}
        \sim
        \pk*{X(a)>u},
        \qquad u\to\infty.
\]
Since   $\lim_{t\to\infty}R(s,t)=2s>0$      this process does not
satisfy the additional condition in
\cite[Proposition~3.1]{chakrabarty2018asymptotic}.

\end{example}

\begin{example} 
\label{34}
Let $X$ be the fractional Brownian sheet with common Hurst index
$H\in(1/2,1)$ and covariance
\[
        R(s,t)
        =
        \prod_{j=1}^2
        \frac{s_j^{2H}+t_j^{2H}-|s_j-t_j|^{2H}}2.
\]
On $\T=[a,b]^2$, where $0<a<b<\infty$, put $t_0=(a,a)$.  The
normalised potential of $\delta_{t_0}$ is
\[
        m(t)=m_H(t_1)m_H(t_2),
        \qquad
        m_H(x):=
        \frac{x^{2H}+a^{2H}-(x-a)^{2H}}{2a^{2H}}.
\]
Strict superadditivity of $x\mapsto x^{2H}$ yields $m_H(x)>1$ for
$x>a$.  Therefore we obtain
\[
        \mu_*=\delta_{t_0},
        \qquad
        S=E=\{t_0\},
        \qquad
        \sigma_*^2=a^{4H}.
\]
Moreover we have
\[
        m_H(a+x)-1=\frac{H}{a}x+o(x),
        \qquad x\downarrow0
\]
and  uniformly on $[a,b]^2$
\[
        \Var(X(t)-X(s))\le C\|t-s\|^{2H}.
\]
The same bound holds for $Z$.   For sufficiently small
$h=(h_1,h_2)$ with $h_1,h_2\ge0$
\[
        m(t_0+h)-1
        \ge \max_{1\le j\le2}\left(m_H(a+h_j)-1\right)
        \ge c\max(h_1,h_2)
        \ge \frac{c}{\sqrt2}\|h\|
\]
for some $c>0$.  Thus $q=1$, $\beta=H$ and $q<2\beta$.  An application of
\Cref{27} yields
\[
        \pk*{\inf_{t\in[a,b]^2}X(t)>u}
        \sim
        \pk*{X(a,a)>u},
        \qquad u\to\infty.
\]
\end{example}

\begin{example} 
\label{35}
Let $\xi_1,\xi_2,\zeta$ be independent standard normal rvs, let
$\T=[-1,1]$, and put
\[
        g(t):=t^2(1-t^2),\qquad
        a_1(t):=\frac{1-t+g(t)}2,\qquad
        a_2(t):=\frac{1+t+g(t)}2.
\]
Setting 
\[
        X(t):=a_1(t)\xi_1+a_2(t)\xi_2+(1-t^2)\zeta,
        \qquad t\in[-1,1]
\]
then for
\[
        \mu_*=\frac12(\delta_{-1}+\delta_1)
\]
we have
\[
        X_{\mu_*}=\frac{\xi_1+\xi_2}{2},
        \qquad
        \sigma_*^2=\frac12,
        \qquad
        m(t)=a_1(t)+a_2(t)=1+g(t).
\]
Hence $\mu_*$ is optimal and
\[
        S=\{-1,1\},\qquad E=\{-1,0,1\},\qquad P=\{0\}.
\]
Any other optimal measure must be supported on $E$.  If its masses at
$-1,0,1$ are $x,y,z$, respectively, then
\[
\begin{aligned}
        \Var\left(\int_\T X d\mu\right)
        &=(x+y/2)^2+(z+y/2)^2+y^2 \ge \frac12+y^2.
\end{aligned}
\]
Thus equality is possible only when $y=0$ and $x=z=1/2$, proving
uniqueness of $\mu_*$.

  Since
\[
        g(-1+h)=2h+O(h^2),
        \qquad
        g(1-h)=2h+O(h^2),
\]
and $Z(t)=(1-t^2)\zeta$, then \[
        \Var(Z(t)-Z(s))
        =(t^2-s^2)^2
        \le4|t-s|^2
\]
and both support contacts are dominated with $q=\beta=1$.  At the passive
contact, $Z(0)=\zeta$ and therefore  
\[
        C=\pk*{Z(0)>0}=\frac12
\] and hence 
\[
        \pk*{\inf_{t\in[-1,1]}X(t)>u}
        \sim
        \frac12 \pk*{\xi_1>u,\ \xi_2>u},
        \qquad u\to\infty.
\]
\end{example}

\subsection{Critical rough and smooth contacts}

We now turn to critical contacts, beginning with a single right-endpoint contact.
Throughout this subsection, $\T=[a,b]$, where $-\infty<a<b<\infty$.
We first assume that $\mu_*=\delta_b$ is optimal, that $R(b,b)>0$, and
\[
        E=S=\{b\}.
\]
Then $\mu_*=\delta_b$ is the unique optimal measure and
\[
        m(t)=\frac{R(t,b)}{R(b,b)},
        \qquad
        m(b)=1,
        \qquad
        m(t)>1,\quad t\in\T\setminus\{b\}.
\]
Since the support is a singleton, there is no distinction between the
optimal-measure residual and the support-regression residual. Thus in the
right-endpoint case
\[
        Z(t):=X(t)-m(t)X(b),
        \qquad t\in\T
\]
is independent of $X(b)$, satisfies $Z(b)=0$ and has continuous sample
paths.\\

We say that $b$ is a tangent-critical endpoint contact with exponent
$\beta\in(0,1]$,
coefficient $\lambda>0$, and tangent process $W$, if the following
conditions hold.  First, for some $C>0$ and $0<r_0\le b-a$
\begin{equation}
\label{36}
        \Var(Z(b-h)-Z(b-h'))
        \le C|h-h'|^{2\beta},
        \qquad 0\le h,h'\le r_0.
\end{equation}
Second, as $h\downarrow0$
\begin{equation}
\label{37}
        m(b-h)-1=\lambda h^{2\beta}+o(h^{2\beta}).
\end{equation}
Third, for every $L>0$ yy the following convergence
of finite-dimensional distributions holds:
\begin{equation}
\label{38}
        \left\{
        r^{-\beta}Z(b-rx),\ 0\le x\le L
        \right\}
        \Rightarrow
        \{W(x),\ 0\le x\le L\}, \quad r\downarrow0,
\end{equation}
where $W=\{W(x), x\ge0\}$ is a centred GP with continuous sample paths and
$W(0)=0$.  When $\beta\in(0,1)$, we call the contact rough-critical.\\

A smooth-critical endpoint contact is the differentiable tangent-critical case
$\beta=1$, in which the residual has a linear $L^2$-tangent.  More precisely,
we assume \eqref{37} with $\beta=1$ and that for some centred Gaussian rv
$\eta$
\[
        \lim_{h\downarrow0}\left\|\frac{Z(b-h)}h-\eta\right\|_{L^2}=0.
\]
We also assume that \eqref{36} holds for some $C>0, 0<r_0\le b-a$ with $\beta=1$.

\begin{theorem}
\label{39}
Under the standing assumptions, let $\T=[a,b]$ with $a<b$.  Assume that
$R(b,b)>0$, that $\mu_*=\delta_b$ is optimal,
and that $E=\{b\}$.  If $b$ is a tangent-critical endpoint contact, then
$M:=\inf_{x\ge0}\left(W(x)+\lambda x^{2\beta}\right)\in(-\infty,0]$ almost surely, and
\[
        \PT
        \sim C_M \pk*{X(b)>u},
        \qquad u\to\infty,
\]
where
\[
        C_M:=\E*{\exp\left(\frac{M}{R(b,b)}\right)}\in(0,1].
\]
\end{theorem}

\begin{corollary}
\label{40}
Under the standing assumptions, suppose that
$R(b,b)>0$, that $\mu_*=\delta_b$ is optimal, and that $E=\{b\}$.
If $b$ is a smooth-critical endpoint contact and
$v^2:=\Var(\eta)$, then
\[
        \PT
        \sim
        \left\{
        \frac12+\frac12
        \left(1+\frac{v^2}{2\lambda R(b,b)}\right)^{-1/2}
        \right\}
        \pk*{X(b)>u},
        \qquad u\to\infty.
\]
\end{corollary}

\begin{remark}
\label{41}
The left-endpoint versions follow by applying the results to
$\widetilde X(t):=X(a+b-t)$.

For an interior contact $t_0$ with $E=S=\{t_0\}$, use
$m(t)=R(t,t_0)/R(t_0,t_0)$ and
$Z(t)=X(t)-m(t)X(t_0)$.  The local increment bound is assumed to hold
across both sides.  Suppose also that for every $L>0$
\[
        \lim_{r\downarrow0}\sup_{|x|\le L}
        \left|
        r^{-2\beta}\{m(t_0+rx)-1\}-d(x)
        \right|=0,
\]
and
\[
        \{r^{-\beta}Z(t_0+rx),\ |x|\le L\}
        \Rightarrow
        \{W(x),\ |x|\le L\}
        \quad\text{in }C[-L,L],
        \qquad r\downarrow0,
\]
and that $m(t_0+h)-1\ge c|h|^{2\beta}$ for some $c>0$ and all sufficiently
small $|h|$.  Then the conclusion of \Cref{39}
holds with $b$ replaced by $t_0$ and with
\[
        M=\inf_{x\in\mathbb R}\left(W(x)+d(x)\right)\in(-\infty,0]
        \quad\text{almost surely}.
\]
In the symmetric smooth case, $W(x)=\eta x$ and
$d(x)=\lambda x^2$, so $M=-\eta^2/(4\lambda)$.  In the asymmetric case
\[
        M=\min\left(M_-,M_+\right),
        \qquad
        M_-:=\inf_{x\le0}\left(W(x)+d(x)\right),
        \qquad
        M_+:=\inf_{x\ge0}\left(W(x)+d(x)\right),
\]
where the two branches come from the joint two-sided tangent.

The interior extension follows from the proof of \Cref{39}, applying
\Cref{56} separately on the two sides and using the assumed joint
convergence on $[-L,L]$.
\end{remark}

\begin{example}
\label{42}
Let $B$ be a standard Brownian motion independent of
$\xi\sim N(0,1)$, and define
\[
        X(t):=\xi(1+\lambda t)+\sigma B(t),
        \qquad t\in[0,1],
        \qquad \lambda,\sigma>0.
\]
For this setting we have 
\[
        \mu_*=\delta_0,
        \qquad S=E=\{0\},
        \qquad m(t)=1+\lambda t,
        \qquad Z(t)=\sigma B(t)
\]
and  for every $L>0$ and all sufficiently small $r>0$
\[
        \{r^{-1/2}Z(rx),\ 0\le x\le L\}
        \stackrel{d}{=}
        \{\sigma B(x),\ 0\le x\le L\}.
\]
Consequently,  $0$ is a rough-critical left-endpoint contact with
$\beta=1/2$ and
\[
        M=\inf_{x\ge0}\left(\sigma B(x)+\lambda x\right).
\]
Since $-M$ is exponential with rate $2\lambda/\sigma^2$ the
left-endpoint version of \Cref{39} yields
\[
        \pk*{\inf_{0\le t\le1}X(t)>u}
        \sim
        \frac{2\lambda}{2\lambda+\sigma^2}
        \pk*{\xi>u},
        \qquad u\to\infty.
\]
\end{example}

\begin{example}
\label{43}  Let $B$ be a standard Brownian motion. Fix
$0<c<1$, $\alpha>1$ and set
\[
        X(t):=\frac{B(t)}{t^\alpha},
        \qquad t\in\T=[c,1].
\]
Since $X(1)=B(1)$ and
\[
        m(t)=\frac{\Cov(X(t),X(1))}{\Var X(1)}
        =t^{1-\alpha},
\]
the optimality criterion yields
\[
        \mu_*=\delta_1,
        \qquad E=S=\{1\},
        \qquad
        Z(t)=t^{-\alpha}(B(t)-tB(1)).
\]
As $h\downarrow0$
\[
        m(1-h)-1=(\alpha-1)h+o(h),
\]
locally uniformly after replacing $h$ by $rx$ with $x$ in a compact
interval.  Moreover, for $s\le t$ near $1$
\[
\begin{aligned}
        \Var(X(t)-X(s))
        &=\frac{t-s}{t^{2\alpha}}
          +s\bigl(t^{-\alpha}-s^{-\alpha}\bigr)^2  \le C|t-s|.
\end{aligned}
\]
The same bound holds for the regression residual $Z$.

For fixed $L<\infty$ and sufficiently small $r$, set
\[
        \widetilde B_r(x)
        :=r^{-1/2}(B(1)-B(1-rx)),
        \qquad 0\le x\le L.
\]
Then $\widetilde B_r$ is a standard Brownian motion on $[0,L]$ and
\[
        r^{-1/2}Z(1-rx)
        =(1-rx)^{-\alpha}
        (-\widetilde B_r(x)+r^{1/2}xB(1))
        \Rightarrow-\widetilde B(x)
\]
in $C[0,L]$, where $\widetilde B$ is a standard Brownian motion.
Hence $1$ is rough-critical with
\[
        \beta=\frac12,
        \qquad \lambda=\alpha-1.
\]
By Brownian symmetry and the calculation in \Cref{42} with
$\sigma=1$ and $\lambda=\alpha-1$, \Cref{39} yields
\[
        \pk*{
        \inf_{t\in[c,1]}\frac{B(t)}{t^\alpha}>u
        }
        \sim
        \frac{2(\alpha-1)}{2\alpha-1}
        \pk*{B(1)>u},
        \qquad u\to\infty.
\]
\end{example}

\begin{remark}
Let $Y\sim N(0,\tau^2)$, $\tau>0$, and let $V$ be a centred
Gaussian random element of $C^1[0,1]$, independent of $Y$, with
$V(1)=0$ almost surely.  If    $m\in C[0,1]$ satisfies
\[
        m(1)=1,
        \qquad
        m(t)>1,\quad 0\le t<1,
        \qquad
        m(1-h)=1+\lambda h^2+o(h^2)
\]
for some $\lambda>0$,    then we have (set  $X(t)=m(t)Y+V(t)$) 
\[
        \mu_*=\delta_1,
        \qquad E=S=\{1\},
        \qquad \sigma_*^2=\tau^2,
        \qquad Z=V
\]
and moreover, Fernique's theorem and the fundamental theorem of calculus yield
\[
        \E*{(V(t)-V(s))^2}
        \le |t-s|^2\E*{\|V'\|_\infty^2},
        \qquad
        \lim_{h\downarrow0}\left\|\frac{V(1-h)}h+V'(1)\right\|_{L^2}=0
\]
and thus \Cref{40} applies with $R(1,1)=\tau^2$ and
$v^2=\Var(V'(1))$.
\end{remark}

\section{Proofs}
\label{44}

\underline{Standing notation for finite-contact proofs}:
Whenever a finite-support optimal measure is under consideration, we write
\[
        A_u:=\{X(t_i)>u,\ i=1,\ldots,k\},
        \qquad
        B_u:=\left\{\inf_{t\in\T}X(t)>u\right\},
\]
and
\[
        \mathbb P_u(\cdot):=\mathbb P(\cdot\mid A_u),
        \qquad
        \bm S_u:=u(\bm X_S-u\bm1).
\]
Since $B_u\subset A_u$, the finite-contact asymptotic proofs below reduce to proving
\[
        \lim_{u\to\infty}\mathbb P_u(B_u)=C
\]
for the appropriate constant $C$.  The residual process $Z$ is independent
of $\bm X_S$ and remains so after conditioning on $A_u$. Setting
\[
        \bm\theta:=\Sigma_*^{-1}\bm1
        =\frac{\bm p_*}{\sigma_*^2}
\]
we have  $\theta_i>0$ for every $i$.  For
$\bm x\in(0,\infty)^k$ a change of variables in the Gaussian density shows
that the conditional density of $\bm S_u$ under $\mathbb P_u$ is
\[
        f_u(\bm x)
        =
        \frac{g_u(\bm x)}
        {\displaystyle\int_{(0,\infty)^k}g_u(\bm y) d\bm y},
        \qquad
        g_u(\bm x)
        :=
        \exp\left(
        -\bm\theta^\top\bm x
        -\frac{1}{2u^2}\bm x^\top\Sigma_*^{-1}\bm x
        \right).
\]
Since $\Sigma_*^{-1}$ is positive definite
\[
        0\le g_u(\bm x)\le \exp\left(-\bm\theta^\top\bm x\right),
\]
and the function on the right is integrable on $(0,\infty)^k$.  Dominated
convergence therefore yields under $\mathbb P_u$ the following weak convergence 
\begin{equation}
\label{45}
        \bm S_u\Rightarrow
        \bm{\mathcal E}:=(\mathcal E_1,\ldots,\mathcal E_k),
        \qquad u\to\infty,
\end{equation}
where the coordinates of $\bm{\mathcal E}$ are independent exponential
random variables with respective rates $\theta_1,\ldots,\theta_k$.  In
particular, $\bm S_u$ is tight under $\mathbb P_u$.

In the one-point endpoint proofs below we take $t_1=b$.  Then
$A_u=\{X(b)>u\}$, $S_u=u(X(b)-u)$, and the preceding convergence holds
with exponential rate
\[
        \theta=\frac1{R(b,b)}>0.
\]

\subsection{Proofs for full contact}
\label{46}

\begin{proof}[Proof of \Cref{18}]
Set
\[
        V:=-\inf_{t\in\T}Z_*(t),
        \qquad
        r(\eps):=\pk*{V<\eps},
        \qquad \eps>0,
\]
and write $s:=\sigma_*^2$.  Path continuity and compactness of $\T$ make
$V$ finite almost surely, while \eqref{11} yields
$V\ge0$.  Since $E=\T$
\[
        \IX=X_{\mu_*}-V,
\]
where $X_{\mu_*}\sim N(0,s)$ is independent of $V$.  Consequently, we have
\[
        f_\T(u)
        =
        \frac1{\sqrt{2\pi s}}
        \E*{
        \exp\left(-\frac{(u+V)^2}{2s}\right)},
\]
with $f_\T$ being a strict positive  density of $X_{\mu_*}$.   Every derivative of the Gaussian density
is bounded and hence   differentiation under the expectation also shows that
$f_\T\in C^\infty(\mathbb R)$. Setting 
\[
        I(u)
        :=
        \int_0^\infty
        \exp\left(
        -\frac{u\eps}{s}
        -\frac{\eps^2}{2s}
        \right)
        r(\eps) d\eps
\]
and by  \eqref{12} 
\begin{equation}
\label{47}
        p_\T(u)
        =
        \frac{\exp\left(-u^2/(2s)\right)}{\sqrt{2\pi s}}I(u).
\end{equation}

Since zero belongs to the support of the centred Gaussian law of
$Z_*$ on $C(\T)$, we have $r(\delta)>0$ for every $\delta>0$.  Define the finite measure
\[
        \nu(d\eps):=\exp\left(-\eps^2/(2s)\right)r(\eps) d\eps,
        \qquad \eps>0,
\]
and the probability measure
\[
        \Pi_u(d\eps):=
        \frac{\exp\left(-u\eps/s\right)\nu(d\eps)}{I(u)}.
\]
For every fixed $\delta>0$, we have $\nu((0,\delta/2))>0$ and
\[
        I(u)\ge \exp\left(-u\delta/(2s)\right)\nu((0,\delta/2)),
        \qquad u>0.
\]
Since $0\le r\le1$, the second moment of $\nu$ is finite.  Splitting at
$\delta$ we therefore obtain
\[
\begin{aligned}
        \int_0^\infty\eps^2 \Pi_u(d\eps)
        &\le \delta^2+
        \frac{\exp\left(-u\delta/s\right)\int_0^\infty\eps^2 \nu(d\eps)}
             {\exp\left(-u\delta/(2s)\right)\nu((0,\delta/2))}\\
        &\le \delta^2+C_\delta \exp\left(-u\delta/(2s)\right).
\end{aligned}
\]
Letting first $u\to\infty$ and then $\delta\downarrow0$ yields
\begin{equation}
\label{48}
        \lim_{u\to\infty}\int_0^\infty\eps^2 \Pi_u(d\eps)=0.
\end{equation}

Write $\bar\eps_u:=\mathbb E_{\Pi_u}\eps$ and
$\tau_u^2:=\operatorname{Var}_{\Pi_u}(\eps)$.  By
\eqref{48}, both tend to zero.  Differentiation
under the integral sign is justified by the Gaussian factor and $0\le r\le1$,
and yields
\[
 (\log I)'(u)=-\frac{\bar\eps_u}{s},\qquad
 (\log I)''(u)=\frac{\tau_u^2}{s^2}.
\]
Since $h_\T=-(\log p_\T)'$,
\eqref{47} therefore yields
\[
 h_\T(u)=\frac{u+\bar\eps_u}{s},\qquad
 h_\T'(u)=\frac1s-\frac{\tau_u^2}{s^2}.
\]
These identities prove \eqref{19}.  Finally, we obtain 
\[
 a_\T(u)=\frac{s}{u+\bar\eps_u}
        =\frac{s}{u}+o(u^{-2}),\qquad
 a_\T'(u)=-\frac{s-\tau_u^2}{(u+\bar\eps_u)^2}
        \sim-\frac{s}{u^2}
\]
establishing    \eqref{20}.
\end{proof}

\begin{lemma}
\label{49}
Assume that $E=\T$ and $\sigma_*^2>0$.
If  for some $\beta,A_1,A_2,\eps_0>0$
\begin{equation}
\label{50}
        \exp\left(-A_2\eps^{-1/\beta}\right)
        \le
        \pk*{\IZ>-\eps}
        \le
        \exp\left(-A_1\eps^{-1/\beta}\right),
        \qquad 0<\eps<\eps_0,
\end{equation}
then there exist constants $0<C_-\le C_+<\infty$ such that for all
sufficiently large $u$
\begin{equation}
\label{51}
        -C_+u^{1/(\beta+1)}
        \le
        \log \PT+\frac{u^2}{2\sigma_*^2}
        \le
        -C_-u^{1/(\beta+1)}.
\end{equation}
\end{lemma}

\begin{proof}[Proof of \Cref{49}]
Let $I_u$ denote the integral in \eqref{12} and put $\kappa:=1/(\beta+1)$ and
$\eps_u:=u^{-\beta/(\beta+1)}$.  For all large $u$,
$2\eps_u<\eps_0$.  Splitting at $\eps_u$, using monotonicity of $r$
on the first part and $r\le1$ on the second, implies
\[
\begin{aligned}
        I_u
        &\le \frac{\sigma_*^2}{u}
        \left[r(\eps_u)+\exp\left(-u\eps_u/\sigma_*^2\right)\right]\le \frac{\sigma_*^2}{u}
        \left[\exp\left(-A_1u^\kappa\right)+\exp\left(-u^\kappa/\sigma_*^2\right)\right].
\end{aligned}
\]
Hence $\log I_u\le-c_2u^\kappa$ for some $c_2>0$ and all large $u$.

For the lower bound monotonicity and
\eqref{50} yield
$r(\eps)\ge r(\eps_u)\ge \exp\left(-A_2u^\kappa\right)$ on
$[\eps_u,2\eps_u]$.  Consequently, we obtain
\[
        I_u
        \ge
        \eps_u
        \exp\left(-\left(A_2+\frac{2}{\sigma_*^2}\right)u^\kappa-1\right)
\]
for all sufficiently large $u$; here the final constant bounds
$\eps^2/(2\sigma_*^2)$ uniformly on
$[\eps_u,2\eps_u]$.  Since
$|\log\eps_u|=O(\log u)=o(u^\kappa)$,
we obtain $\log I_u\ge-c_3u^\kappa$.  Finally
\eqref{12} implies
\[
        \log\PT+\frac{u^2}{2\sigma_*^2}
        =\log I_u-\log(\sqrt{2\pi}\sigma_*),
\]
and the claimed   bound follows after adjusting the constants.
\end{proof}

\begin{proof}[Proof of \Cref{21}]
Let $\mu_*$ and $c_\mu>0$ be as specified before \Cref{21}, and put
$Y:=X_{\mu_*}$. We utilise next the decomposition \eqref{10} and 
identity \eqref{11}.
  For some
$c,C,\varepsilon_0>0$, and all $0<\varepsilon<\varepsilon_0$
\begin{equation}
\label{52}
\begin{split}
 \pk*{\sup_{0\le t\le b}|X(t)|\le\varepsilon}
 &\ge e^{ -C\varepsilon^{-1/H}},\quad 
 \pk*{\int_a^b|X(t)| dt\le\varepsilon}
 \le e^{-c\varepsilon^{-1/H}}.
\end{split}
\end{equation}
The first inequality is the sup-norm small-deviation estimate of
\cite{LiLinde1998}; the second follows from
\cite[Theorem~5.7 and equation~(5.13), with $N=D=q=1$]{LifshitsLindeShi2006}.
For the latter application, normalised Lebesgue measure on $[a,b]$ is
the self-similar measure generated by the two half-interval similarities;
the strong open set condition holds with open set $(a,b)$.  The result therefore
applies directly to this interval also when $a>0$ and normalisation
only changes the constants.

Put $\|f\|_1:=\int_a^b|f(t)| dt$.
If $\sup_{0\le t\le b}|X(t)|\le\varepsilon/3$, then for
$t\in[a,b]$
\[
        |Z_*(t)|
        =
        \left|\int_a^b\{X(t)-X(s)\} \mu_*(ds)\right|
        \le\frac{2\varepsilon}{3}.
\]
The first estimate in \eqref{52} therefore implies that, for some
$0<A_2<\infty$ and all sufficiently small $\varepsilon>0$,
\[
        r(\varepsilon)
        \ge
        \pk*{\sup_{0\le t\le b}|X(t)|
                  \le\frac{\varepsilon}{3}}
        \ge
        \exp\left(-A_2\varepsilon^{-1/H}\right).
\]

For the upper bound on $r$, the anchoring identity implies that on
$\{\IZ>-\varepsilon\}$
\[
 c_\mu\|Z_*\|_1
 \le\int_a^b|Z_*(t)| \mu_*(dt)
 =2\int_a^b(Z_*(t))_- \mu_*(dt)
 \le2\varepsilon.
\]
  Intersecting this event with
$\{|Y|\le\varepsilon\}$, and using $X=Y+Z_*$ and independence yields (set $K:=2/c_\mu+b-a$)
\[
 r(\varepsilon) \pk*{|Y|\le\varepsilon}
 \le\pk*{\|X\|_1\le K\varepsilon}.
\]
Since $\pk*{|Y|\le\varepsilon}\asymp\varepsilon$, the second
estimate in \eqref{52} implies
\[
 r(\varepsilon)
 \le C_1\varepsilon^{-1}
       \exp\left(-c_1\varepsilon^{-1/H}\right)
 \le\exp\left(-A_1\varepsilon^{-1/H}\right)
\]
for some $A_1>0$ and all sufficiently small $\varepsilon>0$.
Thus both inequalities in \eqref{50}
hold with $\beta=H$ and \Cref{49} establishes the claim.  
\end{proof}

\subsection{Proofs for finite contact}
\label{53}

\begin{lemma}
\label{54}
Assume that
\[
        \mu_*=\sum_{i=1}^k p_{*,i}\delta_{t_i},
        \qquad p_{*,i}>0,\quad i=1,\ldots,k,
\]
is the unique optimal probability measure.  If
$P=E\setminus S=\{s_1,\ldots,s_\ell\}$ is nonempty, then
\[
        \pk*{Z(s_j)\ge0,\ j=1,\ldots,\ell}
        =
        \pk*{Z(s_j)>0,\ j=1,\ldots,\ell}
        >0.
\]
\end{lemma}

\begin{proof}
As shown at the start of \Cref{24}, uniqueness implies that $\Sigma_*$
is nonsingular, so the regression representation \eqref{16} applies.  Define
\[
        \bm\Delta:=(Z(s_1),\ldots,Z(s_\ell))^\top
\]
and let $L\subseteq\mathbb R^\ell$ be the linear support of its centred
Gaussian law.  Suppose that
\[
        \pk*{\Delta_j>0,\ j=1,\ldots,\ell}=0.
\]
Since the Gaussian law has full support on $L$,
$L\cap(0,\infty)^\ell=\varnothing$.  Gordan's alternative provides
a nonzero vector $\bm\lambda\in[0,\infty)^\ell$ such that
\[
        \bm\lambda^\top\bm x=0,
        \qquad \bm x\in L.
\]
Consequently, $\bm\lambda^\top\bm\Delta=0$ almost surely.

Let
\[
        \bm A
        :=
        (a_i(s_j))_{\substack{1\le i\le k\\1\le j\le\ell}},
        \qquad
        \bm c:=\bm A\bm\lambda.
\]
Since $m(s_j)=\sum_{i=1}^k a_i(s_j)=1$
\[
        \sum_{i=1}^k c_i
        =
        \sum_{j=1}^\ell\lambda_jm(s_j)
        =
        \sum_{j=1}^\ell\lambda_j.
\]
As every $p_{*,i}>0$, for all sufficiently small $\varepsilon>0$
\[
        \mu_\varepsilon
        :=
        \sum_{i=1}^k
        (p_{*,i}-\varepsilon c_i)\delta_{t_i}
        +
        \varepsilon\sum_{j=1}^\ell\lambda_j\delta_{s_j}
\]
is a probability measure.    Writing
\[
        \bm X_P:=(X(s_1),\ldots,X(s_\ell))^\top
\]
we have $\bm X_P=\bm A^\top\bm X_S+\bm\Delta$ and hence
\[
\begin{aligned}
        \int_\T X(t) \mu_\varepsilon(dt)
        -
        \int_\T X(t) \mu_*(dt)
        &=
        \varepsilon
        \left(
        \bm\lambda^\top\bm X_P-\bm c^\top\bm X_S
        \right)                                                                  \\
        &=
        \varepsilon\bm\lambda^\top
        \left(\bm X_P-\bm A^\top\bm X_S\right)                     \\
        &=
        \varepsilon\bm\lambda^\top\bm\Delta
        =0
        \qquad\text{a.s.}
\end{aligned}
\]
Thus $\mu_\varepsilon$ is another optimal measure.  It is distinct from
$\mu_*$, because $\bm\lambda\ne0$ and $s_j\notin S$, contradicting
uniqueness.  Therefore we have 
\[
        \pk*{\Delta_j>0,\ j=1,\ldots,\ell}>0.
\]

In particular, no $\Delta_j$ is identically zero.  Each coordinate is
therefore a nondegenerate Gaussian rv and hence  
\[
        \pk*{\Delta_j=0}=0,
        \qquad j=1,\ldots,\ell
\] establishing the proof. 
\end{proof}

\medskip

\begin{proof}[Proof of \Cref{27}]
\Cref{54} yields $C>0$ and almost surely no passive value is zero.  When $P=\varnothing$, these
conditions are void and $C=1$.

For each $i$, choose an even integer $N_i\ge2$ and an exponent
$\gamma_i$ such that
\[
        N_i\left(\beta_i-\frac{q_i}{2}\right)>d,
        \qquad
        \frac{q_i}{2}<\gamma_i<
        \min\left(q_i,\beta_i-\frac{d}{N_i}\right).
\]
On sufficiently small relative neighbourhoods $U_i$ of $t_i$ where
\eqref{26} holds, Gaussian moment estimates yield
\[
        \E*{|Z(t_i+h)-Z(t_i+h')|^{N_i}}
        \le C_{i,N_i}\|h-h'\|^{N_i\beta_i}.
\]
On a smaller compact relative neighbourhood of $t_i$, choose, for each $n$,
one representative from each occupied dyadic cube of side $2^{-n}$, taking
$t_i$ as the representative of its cube.
Joining representatives in neighbouring cubes at scale $n$ and in
parent--child cubes at scales $n$ and $n+1$
gives a set $\mathcal E_n$ of $O(2^{nd})$ edges, each of length
$O(2^{-n})$.  The moment bound and a union bound yield
\[
 \pk*{
 \max_{(s,t)\in\mathcal E_n}|Z(s)-Z(t)|>2^{-n\gamma_i}
 }
 \le C_i'2^{-n\{N_i(\beta_i-\gamma_i)-d\}}.
\]
The right-hand side is summable.  Borel--Cantelli and dyadic chaining,
using the given continuous version and $Z(t_i)=0$, therefore yield,
after shrinking $U_i$, an almost surely finite rv $K_i$ such that
\[
        |Z(t_i+h)|\le K_i\|h\|^{\gamma_i},
        \qquad t_i+h\in U_i.
\]

For $\bm x\in(0,\infty)^k$ and $z\in C(K)$, put
\[
        H_u(\bm x,z):=
        \mathbf1\left\{
        \inf_{t\in K}
        \left[u\{m(t)-1\}+z(t)+\frac{\bm a(t)^\top\bm x}{u}\right]>0
        \right\}.
\]
For every fixed $\bm x\in(0,\infty)^k$, we claim that
\[
        \lim_{u\to\infty}H_u(\bm x,Z)=\mathbf1\{Z(s)>0,\ s\in P\}
        \qquad\text{a.s.}
\]
Fix a residual path satisfying the preceding bounds and having no zero
passive value.  If $Z(s)<0$ at some $s\in P$, evaluating at $s$, where
$m(s)=1$, shows that $H_u(\bm x,Z)=0$ for all large $u$.

Suppose instead that all passive values are positive.  Near each active
contact, continuity and $\bm a(t_i)=\bm e_i$ yield
$\bm a(t)^\top\bm x\ge x_i/2$.  Since $q_i/2<\gamma_i<q_i$
\[
        \sup_{r\ge0}\left(K_i r^{\gamma_i}-\lambda_i u r^{q_i}\right)
        =O\left(u^{-\gamma_i/(q_i-\gamma_i)}\right)=o(u^{-1}).
\]
Thus the expression defining $H_u$ is bounded below by
$x_i/(2u)-o(u^{-1})>0$, uniformly on a sufficiently small relative
neighbourhood of $t_i$.  Near each passive contact, continuity makes $Z$
bounded below by a positive constant; since $m\ge1$ and $\bm a$ is
bounded, the same expression is positive there for all large $u$.
The complement of these relative open neighbourhoods is compact and disjoint
from $E$.  If nonempty, $m-1$ is bounded below by a strictly positive constant there, while $Z$ and $\bm a$
are bounded.  Positivity follows there as well, proving the claim.

Using the conditional density $f_u$ and its unnormalised version $g_u$
from the standing notation, independence of $Z$ and $\bm X_S$ yields
\[
        \mathbb P_u(B_u)
        =
        \frac{\displaystyle
        \int_{(0,\infty)^k}g_u(\bm x)
        \E{ H_u(\bm x,Z)} d\bm x}
        {\displaystyle\int_{(0,\infty)^k}g_u(\bm x) d\bm x}.
\]
The pointwise convergence above and $0\le H_u\le1$ imply
$\lim_{u\to \infty}\E {H_u(\bm x,Z)}= C$.  Since
\[
        \lim_{u\to\infty}g_u(\bm x)=\exp\left(-\bm\theta^\top\bm x\right)
\]
the integrable bound in the standing notation permits dominated convergence,
yielding
$ \lim_{u\to\infty}\mathbb P_u(B_u)= C$.  This proves the first
asymptotic relation.  The same Gaussian change of variables implies
\[
        \mathbb P(A_u)
        =\frac{u^{-k}\exp\left(-u^2/(2\sigma_*^2)\right)}
               {(2\pi)^{k/2}\sqrt{\det\Sigma_*}}
          \int_{(0,\infty)^k}g_u(\bm x) d\bm x,
\]
and dominated convergence sends the integral to
$\prod_{i=1}^k\theta_i^{-1}$, proving the explicit prefactor.

\end{proof}

\medskip

\begin{lemma}
\label{55}
Fix $\beta\in(0,1]$.  Let $D\subset[0,\varrho]$ be a nonempty compact set containing
$0$, where $\varrho\ge0$.  Let $G(x), x\in D$ be a centred GP
with continuous sample paths, $G(0)=0$, and suppose that for some $0\le A<\infty$
\[
  \E*{(G(x)-G(y))^2}
  \le
  A^2|x-y|^{2\beta},
  \qquad x,y\in D.
\]
Then there exists a constant $K_\beta<\infty$, depending only on $\beta$, such
that uniformly in $A$, $\varrho$, and $D$
\[
  \E*{\sup_{x\in D}|G(x)|}
  \le
  K_\beta A \varrho^\beta, \qquad
  \sup_{x\in D}\Var G(x)
  \le
  A^2\varrho^{2\beta}.
\]
\end{lemma}
\begin{proof}
If $A=0$ or $\varrho=0$, the result is immediate from anchoring and
continuity.  Otherwise, the canonical metric satisfies
\[
        d_G(x,y)\le A|x-y|^\beta,
        \qquad
        N(D,d_G,\eps)\le1+\varrho(A/\eps)^{1/\beta}.
\]
For the covering bound, a maximal $(\eps/A)^{1/\beta}$-separated
subset of $D$ in the Euclidean metric is an $\eps$-net in $d_G$ and has at most
$1+\varrho(A/\eps)^{1/\beta}$ points.
Since $G(0)=0$, the variance bound follows by taking $y=0$.
Dudley's entropy bound and Gaussian symmetry yield, with a universal
constant $K$
\[
\begin{aligned}
        \E*{\sup_{x\in D}|G(x)|}
        &\le
        2K\int_0^{A\varrho^\beta}
        \sqrt{\log\left(1+\varrho(A/\eps)^{1/\beta}\right)} d\eps\\
        &=2KA\varrho^\beta\int_0^1
        \sqrt{\log\left(1+s^{-1/\beta}\right)} ds
        \le K_\beta A\varrho^\beta.
\end{aligned}
\]
The integral is finite for every $\beta>0$.
\end{proof}

\medskip

\begin{lemma}
\label{56}
Fix $\beta\in(0,1]$, $A\ge0$, $c>0$, and $R\in(0,\infty]$.
Let $G$ be a centred GP on $[0,R]$, interpreted as $[0,\infty)$
when $R=\infty$, with continuous sample paths, $G(0)=0$, and satisfying
the increment bound of \Cref{55} for all $x,y$ in its domain.
Let $d$ be continuous on the same domain and satisfy
$d(x)\ge cx^{2\beta}$.  There exists a finite constant
$K=K(\beta,A,c)$, independent of $R,G,d$, such that for $0<L\le R$
\[
        \pk*{\inf_{L\le x\le R}\left(G(x)+d(x)\right)\le0}
        \le K L^{-\beta}.
\]
If $R=\infty$, then also
$\inf_{x\ge0}\left(G(x)+d(x)\right)>-\infty$ almost surely.
\end{lemma}
\begin{proof}
Put $Q_n=2^nL$, and retain the intervals
$I_n=[Q_n,2Q_n]\cap[0,R]$ that are nonempty.  Utilising 
\Cref{55} and Markov's inequality we obtain (set $\varrho_n=\min\left(2Q_n,R\right)$)
\[
\begin{aligned}
        \pk*{
        \sup_{0\le x\le\varrho_n}|G(x)|
        \ge\frac c2Q_n^{2\beta}}
        &\le\frac{2K_\beta A\varrho_n^\beta}{cQ_n^{2\beta}} \le\frac{2^{\beta+1}K_\beta A}{c} Q_n^{-\beta}.
\end{aligned}
\]
Outside the union of these events, for every $x\in I_n$
\[
        G(x)+d(x)>\frac c2Q_n^{2\beta}\ge\frac c2L^{2\beta}.
\]
Consequently, summing the geometric series yields the asserted estimate,
with $$K=2^{\beta+1}K_\beta A/[c(1-2^{-\beta})].$$
When $R=\infty$, the events
$\{\inf_{x\ge L}(G(x)+d(x))\le0\}$ decrease as $L\to\infty$, and their
probabilities tend to zero.  Almost surely the process is therefore
positive outside some compact interval.  Continuity on that interval
proves the final assertion.
\end{proof}

\medskip

\begin{proof}[Proof of \Cref{39}]
By the standing reduction, it remains to identify the limit of
$\mathbb P_u(B_u)$.  We use the one-point overshoot convergence stated
above and put
\[
        r_u:=u^{-1/\beta},
        \qquad
        T_u:=\sup_{t\in[a,b]}
        \frac{-uZ(t)-u^2\{m(t)-1\}}{m(t)}.
\]
Since $m\ge1$ and $Z(b)=0$, we have $T_u\ge0$.  Gaussian regression,
compactness and continuity yield the exact identity
\[
        B_u=\{S_u>T_u\}.
\]
The rv $T_u$ depends only on $Z$, so it is independent of $S_u$ and
its law is unchanged under $\mathbb P_u$.

Since all processes in \eqref{38} are centred Gaussian, taking
$r\downarrow0$ in \eqref{36} shows that $W$ satisfies the increment
bound of \Cref{55} on $[0,\infty)$ with $A=\sqrt C$.
\Cref{56}, applied with $G=W$,
$d(x)=\lambda x^{2\beta}$, and $R=\infty$, yields $M>-\infty$
almost surely.  Since $W(0)=0$, also $M\le0$ and hence
$C_M\in(0,1]$.

For every fixed $L<\infty$, the expansion
\eqref{37} yields, uniformly on $0\le x\le L$
\[
        \lim_{u\to\infty}\sup_{0\le x\le L}
        \left|u^2(m(b-r_ux)-1)-\lambda x^{2\beta}\right|=0,
        \qquad
        \lim_{u\to\infty}\sup_{0\le x\le L}|m(b-r_ux)-1|=0.
\]
Indeed, the error in the first convergence is bounded by
\[
        L^{2\beta}\sup_{0<h\le r_uL}
        \left|\frac{m(b-h)-1}{h^{2\beta}}-\lambda\right|.
\]
Moreover, the finite-dimensional convergence
\eqref{38} can be upgraded to weak convergence in
$C[0,L]$.  To see this, choose
an even integer $p$ with $p\beta>1$.  The increment bound and Gaussian
moments imply, uniformly for all sufficiently large $u$
\[
        \E*{|uZ(b-r_ux)-uZ(b-r_uy)|^p}
        \le C_p|x-y|^{p\beta},
        \qquad 0\le x,y\le L.
\]
The Kolmogorov--Chentsov criterion therefore yields tightness. Set
\[
\begin{aligned}
        T_{u,L}
        &:=\sup_{0\le x\le L}
        \frac{-uZ(b-r_ux)-u^2\{m(b-r_ux)-1\}}{m(b-r_ux)},\quad 
        D_L:=\inf_{0\le x\le L}\left(W(x)+\lambda x^{2\beta}\right).
\end{aligned}
\]
The continuous mapping theorem yields $T_{u,L}\Rightarrow-D_L$ as $u\to\infty$.
Indeed, the rescaled residual process, drift, and denominator converge jointly
in $C[0,L]^3$ to $(W,\lambda x^{2\beta},1)$, and the map
$(f,d,m)\mapsto\sup_{[0,L]}\{(-f-d)/m\}$ is continuous whenever $\inf m>0$.

Choose $0<r_0<b-a$ so that the local increment bound holds and
$m(b-h)-1\ge(\lambda/2)h^{2\beta}$ for $0\le h\le r_0$.
By \eqref{36}, on $[0,r_0/r_u]$ the process $G_u(x)=uZ(b-r_ux)$
satisfies the increment bound of \Cref{55} with $A=\sqrt C$, uniformly
in $u$, while
\[
        u^2(m(b-r_ux)-1)\ge\frac\lambda2x^{2\beta}.
\]
Thus \Cref{56} bounds the probability of a nonpositive
value on $Lr_u\le h\le r_0$ by $K L^{-\beta}$, uniformly in all
sufficiently large $u$.  Since $E=\{b\}$
\[
        d_0:=\inf_{a\le t\le b-r_0}\left(m(t)-1\right)>0.
\]
As $T_{u,L}\ge0$, an outside point can change the supremum only if its
numerator is positive.  Hence, for all sufficiently large $u$
\[
\begin{aligned}
        \pk*{T_u\ne T_{u,L}}
        &\le
        \pk*{
        \inf_{Lr_u\le h\le r_0}
        \left[u^2(m(b-h)-1)+uZ(b-h)\right]\le0
        } +
        \pk*{\sup_{t\in[a,b]}|Z(t)|\ge ud_0}.
\end{aligned}
\]
The last probability tends to zero by path continuity and compactness.
Together with \Cref{56} this proves
\[
        \lim_{L\to\infty}\limsup_{u\to\infty}
        \pk*{T_u\ne T_{u,L}}=0.
\]
Since $D_L\downarrow M>-\infty$ we conclude the weak convergence of
$T_u\Rightarrow-M$ as $u\to\infty$.  Independence and the one-point overshoot limit
therefore yield, under $\mathbb P_u$
\[
        (S_u,T_u)\Rightarrow(\mathcal E,-M), \quad u\to\infty,
\]
where $\mathcal E$ is independent of $M$ and is exponential with rate
$1/R(b,b)$.  In particular, $\pk*{\mathcal E=-M}=0$, and hence
\[
        \lim_{u\to\infty}\mathbb P_u(B_u)
        =\pk*{\mathcal E>-M}
        =\E*{\exp\left(\frac{M}{R(b,b)}\right)}.
\]
This proves the assertion.
\end{proof}

\medskip

\begin{proof}[Proof of \Cref{40}]
For every fixed $x\ge0$
\[
        \lim_{u\to\infty}\left\|uZ(b-x/u)-\eta x\right\|_{L^2}=0.
\]
For $x>0$, this follows from
$uZ(b-x/u)=xZ(b-x/u)/(x/u)$ and the assumed $L^2$-tangent;
for $x=0$, both terms vanish.  This yields convergence of the
finite-dimensional distributions.  Together with the assumed local
mean-square Lipschitz bound, this shows that $b$ is a tangent-critical
endpoint contact with
\[
        \beta=1,
        \qquad
        W(x)=\eta x,\quad x\ge0.
\]
Its variational term is
\[
\begin{aligned}
        M_{\rm sm}
        &:=
        \inf_{x\ge0}\left(\eta x+\lambda x^2\right)        =
        -\frac{(\max\left(-\eta,0\right))^2}{4\lambda}.
\end{aligned}
\]
The conclusion of \Cref{39} therefore applies.
For $s:=R(b,b)$ a standard Gaussian calculation yields the following
identity valid also when $v=0$:
\[
\begin{aligned}
\E*{e^{-\frac{(\max\left(-\eta,0\right))^2}{4\lambda s}}}
&=\frac12+\frac12\E{e^ {-\frac{\eta^2}{4\lambda s}}}=\frac12+\frac12
  \left(1+\frac{v^2}{2\lambda s}\right)^{-1/2}
\end{aligned}
\]
establishing the claim. 
\end{proof}

\appendix
\section{Auxiliary Results}
Let $\mathcal T$ be a nonempty compact metric space and let
$X(t), t\in\mathcal T$ be a centred real valued GP with
almost surely continuous sample paths and covariance function $R$. For a finite signed measure $\gamma$ on $\mathcal T$, define the energy
associated with the covariance kernel $R$ by
\[
        \mathcal E_R(\gamma)
        =\iint_{\mathcal T^2}R(s,t) \gamma(ds)\gamma(dt)\ge0.
\]
For $\mu\in\mathcal M_1(\mathcal T)$ define its potential and energy by
\[
        \phi_\mu(t)=\int_{\mathcal T}R(s,t) \mu(ds),
        \qquad
        \sigma^2(\mu):=\mathcal E_R(\mu).
\]
We call $\mu$ an optimal measure if
\[
        \sigma^2(\mu)=
        \inf_{\nu\in\mathcal M_1(\mathcal T)}\sigma^2(\nu).
\]

We use the following standard optimality criterion for minimum-energy
probability measures; see \cite[Theorem~3.1]{hashorvaNovikov2026highminima}.
Related formulations appear in
\cite[Theorem~4.3 and Remark~4.5]{adler2014existence},
\cite[Theorem~2.1]{selk2024large}, and
\cite[Theorem~1.1]{pronzato2021minimum}.

\begin{lemma}
\label{57}
An optimal measure exists.  A probability measure
$\mu\in\mathcal M_1(\mathcal T)$ is optimal if and only if
\[
        \phi_\mu(t)\ge \sigma^2(\mu),
        \qquad t\in\mathcal T.
\]
For every optimal $\mu$, equality holds on $\supp(\mu)$. If further  
$\mathcal E_R(\gamma)>0$ for every nonzero finite signed measure $\gamma$ with
$\gamma(\mathcal T)=0$, then the optimal measure is unique.
\end{lemma}

\begin{corollary}
\label{58}
Let $S=\{t_1,\ldots,t_m\}\subset\mathcal T$, assume that
$\Sigma=(R(t_i,t_j))_{i,j=1}^m$ is invertible and set
\[
        w:=
        \frac{\Sigma^{-1}\mathbf 1}{\mathbf 1^\top \Sigma^{-1}\mathbf 1},
        \qquad
        c:=
        \frac{1}{\mathbf 1^\top \Sigma^{-1}\mathbf 1},
        \qquad
        E:=\left\{t\in\mathcal T:
        \sum_{j=1}^m w_jR(t_j,t)=c\right\}.
\]
If all components of $w$ are positive and
\[
        \sum_{j=1}^m w_jR(t_j,t)\ge c,
        \qquad t\in\mathcal T,
\]
then
$
        \mu_*:=\sum_{j=1}^m w_j\delta_{t_j}
$
is an optimal measure with support exactly $S$. If, in addition, either
$\mathcal E_R(\gamma)> 0$ for every nonzero finite signed measure $\gamma$ on $\mathcal T$ with
$\gamma(\mathcal T)=0$, or $E=S$, then
$\mu_*$ is the unique optimal measure.
Conversely, if an optimal measure has support   $S$, then it equals
$\sum_{j=1}^m w_j\delta_{t_j}$ and $w$ has
positive components.
\end{corollary}

\begin{proof}
Because $\Sigma$ is an invertible covariance matrix, it is positive
definite. Hence $c>0$, $\mathbf1^\top w=1$,
$\Sigma w=c\mathbf1$, and $\sigma^2(\mu_*)=c$.  The assumed potential
inequality and \Cref{57} yield optimality; the
strict-energy alternative yields uniqueness.

Suppose instead that $E=S$, and let $\nu$ be another optimal measure.
Since $\phi_{\mu_*}\ge c$ and both measures have energy $c$
\[
0\le\mathcal E_R(\nu-\mu_*)
=2c-2\int_{\mathcal T}\phi_{\mu_*} d\nu
\le0
\]
and hence  $\int_{\mathcal T}(\phi_{\mu_*}-c) d\nu=0$.  Continuity and nonnegativity imply
$\supp(\nu)\subseteq E=S$.  Writing
$\nu=\sum_{j=1}^m v_j\delta_{t_j}$, we obtain
\[
        0=\mathcal E_R(\nu-\mu_*)=(v-w)^\top\Sigma(v-w)
\]
and positive definiteness yields $v=w$.

For the converse, write the optimal measure as
$\mu=\sum_{j=1}^m v_j\delta_{t_j}$ and put $d:=\sigma^2(\mu)$.  Equality of its
potential and energy on its support yields $\Sigma v=d\mathbf1$.
Together with $\mathbf1^\top v=1$, this yields $d=c$ and $v=w$.
Exact support yields $v_j>0$ for every $j$.
\end{proof}

\bibliographystyle{ieeetr}

\bibliography{GG_V47}

\end{document}